\documentclass[11pt, reqno, psamsfonts]{amsart}
\pdfoutput=1

\usepackage{amssymb}
\usepackage{amsthm}
\usepackage{amsmath}
\usepackage{latexsym}
\usepackage[T1]{fontenc}
\usepackage[utf8]{inputenc}
\usepackage[russian, french, english]{babel}

\usepackage{graphicx}
\usepackage{wrapfig}
\usepackage[justification=centering, labelfont=bf]{caption}
\usepackage{mathtools}
\usepackage{tikz}
\usepackage{amsbsy}
\usepackage[inline]{enumitem}
\usepackage{mathrsfs}
\usetikzlibrary{shapes,snakes}
\usetikzlibrary{arrows.meta}
\usetikzlibrary{decorations.pathmorphing}
\usetikzlibrary{patterns}
\usepackage{float}
\usepackage{array}
\usepackage{multicol}
\usepackage{stmaryrd}
\usepackage{cancel} 
\usepackage{lmodern}
\usepackage{mathabx}
\usepackage{upgreek}
\usepackage{titlesec}
\usepackage[spacing=true,kerning=true,babel=true,tracking=true]{microtype}
\usepackage[shortcuts]{extdash}
\usepackage[foot]{amsaddr}
\usepackage[left=1in,right=1in,top=1in,bottom=1in,bindingoffset=0cm]{geometry}
\usepackage{bm}
\usepackage{centernot}
\usepackage{mdframed}
\usepackage[hidelinks]{hyperref}
\usepackage{xspace}
\usepackage{eucal,eufrak}

\usepackage[backend=biber, style=alphabetic, sorting=nyt, maxnames=100,backref=true]{biblatex}
\title{\sffamily Borel fractional colorings of Schreier graphs}
\date{}

\author{Anton~Bernshteyn}
\address{\normalfont School of Mathematics, Georgia Institute of Technology, Atlanta, GA, USA}
\email{bahtoh@gatech.edu}

\thanks{This research is partially supported by the NSF grant DMS-2045412.}

\newtheoremstyle{bfnote}%
{}{}%
{\slshape}{}%
{\bfseries}{\bfseries.}%
{ }%
{\thmname{#1}\thmnumber{ #2}\thmnote{ \ep{\normalfont{}#3}}}

\theoremstyle{bfnote}
\newtheorem{theo}{Theorem}[section]
\newtheorem*{theo*}{Theorem}

\newtheorem{lemma}[theo]{Lemma}

\newtheorem{corl}[theo]{Corollary}

\newtheorem*{corl*}{Corollary}

\theoremstyle{definition}

\newtheorem*{defn*}{Definition}

\newtheorem{ques}[theo]{Question}

\newtheorem*{exmp*}{Example}

\theoremstyle{remark}
\newtheorem*{ques*}{Question}
\newtheorem*{remk*}{Remark}

\newcommand{\0}{\varnothing}
\newcommand{\set}[1]{\{#1\}}
\newcommand{\N}{{\mathbb{N}}}

\newcommand{\F}{\mathbb{F}}

\renewcommand{\epsilon}{\varepsilon}

\renewcommand{\phi}{\varphi}
\renewcommand{\theta}{\vartheta}
\renewcommand{\leq}{\leqslant}
\renewcommand{\geq}{\geqslant}
\newcommand{\defeq}{\coloneqq}

\newcommand{\bemph}[1]{{\normalfont#1}} 
\newcommand{\ep}[1]{\bemph{(}#1\bemph{)}} 

\newcommand{\emphd}[1]{{\fontseries{b}\selectfont\textsf{#1}}}
\newcommand{\acts}{\mathrel{\reflectbox{$\righttoleftarrow$}}}
\newcommand{\racts}{\mathrel{\righttoleftarrow}}
\renewcommand{\G}{\Gamma}
\newcommand{\pto}{\dashrightarrow}
\newcommand{\dom}{\mathrm{dom}}
\newcommand{\Free}{\mathrm{Free}}

\newcommand{\symdif}{\vartriangle}
\newcommand{\rest}[2]{{{#1}\vert_{#2}}}

\numberwithin{equation}{section}

\newenvironment{scproof}[1][]{\begin{proof}[\textsc{\upshape{Proof}}#1]}{\end{proof}}

\titleformat{\section}[block]{\large\bfseries\sffamily}{\thesection.}{1ex}{}
\titleformat{\subsection}[block]{\bfseries\sffamily}{\thesubsection.}{1ex}{}
\titleformat{\subsubsection}[runin]{\itshape}{\bfseries\upshape\thesubsubsection.}{1ex}{}[.---]

\titlespacing*{\section}{0pt}{*3}{*1}
\titlespacing*{\subsection}{0pt}{*3}{*1}
\titlespacing*{\subsubsection}{0pt}{*1.5}{*0}

\renewbibmacro{in:}{}

\renewbibmacro*{volume+number+eid}{%
	\printfield{volume}%
	\setunit*{\addnbspace}
	\printfield{number}%
	\setunit{\addcomma\space}%
	\printfield{eid}}

\DeclareFieldFormat[article]{volume}{\textbf{#1}\space}
\DeclareFieldFormat[article]{number}{\mkbibparens{#1}}

\DeclareFieldFormat{journaltitle}{#1,}
\DeclareFieldFormat[thesis]{title}{\mkbibemph{#1}\addperiod}
\DeclareFieldFormat[article, unpublished, thesis]{title}{\mkbibemph{#1},}
\DeclareFieldFormat[book]{title}{\mkbibemph{#1}\addperiod}
\DeclareFieldFormat[unpublished]{howpublished}{#1, }

\DeclareFieldFormat{pages}{#1}

\DeclareFieldFormat[article]{series}{Ser.~#1\addcomma}

\setlist{topsep=3pt,itemsep=3pt}

\begin{document}


    \maketitle
    
    
    \begin{abstract}
        Let $\G$ be a countable group and let $G$ be the Schreier graph of the free part of the Bernoulli shift $\G \acts 2^\G$ (with respect to some finite subset $F \subseteq \G$). We show that the Borel fractional chromatic number of $G$ is equal to $1$ over the measurable independence number of $G$. As a consequence, we asymptotically determine the Borel fractional chromatic number of $G$ when $\G$ is the free group, answering a question of Meehan.
    \end{abstract}
    
    \section{Definitions and results}
    
    All graphs in this paper are undirected and simple. Recall that for a graph $G$, a subset $I \subseteq V(G)$ is \emphd{$G$-independent} if no two vertices in $I$ are adjacent in $G$. The \emphd{chromatic number} of $G$, denoted by $\chi(G)$, is the least $\ell \in \N$ such that there exist $G$-independent sets $I_1$, \ldots, $I_\ell$ whose union is $V(G)$. \ep{If no such $\ell$ exists, we set $\chi(G) \defeq \infty$.} The sequence $I_1$, \ldots, $I_\ell$ is called an \emphd{$\ell$-coloring} of $G$, where we think of the vertices in $I_i$ as being assigned the color $i$.
    
    Fractional coloring is a well-studied relaxation of graph coloring. For an introduction to this topic, see the book \cite{FracBook} by Scheinerman and Ullman. Given $k \in \N$, the \emphd{$k$-fold chromatic number} of $G$, denoted by $\chi^k(G)$, is the least $\ell \in \N$ such that there are $G$-independent sets $I_1$, \ldots, $I_\ell$ which cover every vertex of $G$ at least $k$ times \ep{such a sequence $I_1$, \ldots, $I_\ell$ is called a \emphd{$k$-fold $\ell$-coloring}}. Note that the sets $I_1$, \ldots, $I_\ell$ need not be distinct. In particular, if $I_1$, \ldots, $I_{\chi(G)}$ is a $\chi(G)$-coloring of $G$, then, by repeating each set $k$ times, we obtain a $k$-fold $k\chi(G)$-coloring, which shows that
    \[
        \chi^k(G) \,\leq\, k \chi(G) \quad \text{for all } k.
    \]
    This inequality can be strict; for example, the $5$-cycle $C_5$ satisfies $\chi(C_5) = 3$ but $\chi^2(C_5) = 5$. It is therefore natural to define the \emphd{fractional chromatic number} $\chi^\ast(G)$ of $G$ by the formula
	\[
	    \chi^\ast(G) \,\defeq\, \inf_{k \geq 1} \frac{\chi^k(G)}{k}.
	\]


	
	In this note we investigate fractional colorings from the standpoint of Borel combinatorics. For a general overview of Borel combinatorics, see the surveys \cite{KechrisMarks} by Kechris and Marks and \cite{Pikh_survey} by Pikhurko. The study of fractional colorings in this setting was initiated by Meehan \cite{Mee}; see also \cite[\S8.6]{KechrisMarks}. 
	We say that a graph $G$ is \emphd{Borel} if $V(G)$ is a standard Borel space and $E(G)$ is a Borel subset of $V(G) \times V(G)$. The \emphd{Borel chromatic number} $\chi_\mathrm{B}(G)$ of $G$ is the least $\ell \in \N$ such that there exist {Borel} $G$-independent sets $I_1$, \ldots, $I_\ell$ whose union is $V(G)$. The \emphd{Borel $k$-fold chromatic number} $\chi^k_{\mathrm{B}}(G)$ is defined analogously,  
	and the \emphd{Borel fractional chromatic number} $\chi^\ast_{\mathrm{B}}(G)$ is 
	\[
	    \chi^\ast_{\mathrm{B}}(G) \,\defeq\, \inf_{k \geq 1} \frac{\chi^k_{\mathrm{B}}(G)}{k}.
	\]
	
	A particularly important class of Borel graphs are Schreier graphs of group actions. Let $\G$ be a countable group with identity element $\mathbf{1}$ and let $F \subseteq \G$ be a finite subset. The \emphd{Cayley graph} $G(\G, F)$ of $\G$ is the graph with vertex set $\G$ in which two distinct group elements $\gamma$, $\delta$ are adjacent if and only if $\gamma =\sigma \delta$ for some $\sigma \in F \cup F^{-1}$. This definition can be extended as follows. Let $\G \acts X$ be a Borel action of $\G$ on a standard Borel space $X$. The action $\G \acts X$ is \emphd{free} if \[\gamma \cdot x \,\neq\, x \quad \text{for all } x \in X \text{ and } \mathbf{1} \neq \gamma \in \G.\] The \emphd{Schreier graph} $G(X,F)$ of an action $\G \acts X$ is the graph with vertex set $X$ in which two distinct points $x$, $y \in X$ are adjacent if and only if $y = \sigma \cdot x$ for some $\sigma \in F \cup F^{-1}$. 
	Note that the Cayley graph $G(\G, F)$ is a special case of this construction corresponding to the left multiplication action $\G \acts \G$. More generally, when the action $\G \acts X$ is free, $G(X, F)$ is obtained by putting a copy of the Cayley graph $G(\G, F)$ onto each orbit.
	
	A crucial example of a Borel action is the \emphd{\ep{Bernoulli} shift} $\G \acts 2^\G$, given by the formula
	\[
	    (\gamma \cdot x)(\delta) \,\defeq\, x(\delta \gamma) \quad \text{for all } x \colon \G \to 2 \text{ and } \gamma,\ \delta \in \G.
	\]
	We use $\beta$ to denote the ``coin flip'' probability measure on $2^\G$, obtained as the product of countably many copies of the uniform probability measure on $2 = \set{0,1}$. Note that $\beta$ is invariant under the shift action. The \emphd{free part} of $2^\G$, denoted by $\Free(2^\G)$, is the set of all points $x \in 2^\G$ with trivial stabilizer. In other words, $\Free(2^\G)$ is the largest subspace of $2^\G$ on which the shift action is free. It is easy to see that the shift action $\G \acts 2^\G$ is free $\beta$-almost everywhere, i.e., $\beta(\Free(2^\G)) = 1$.
	
	
	
	Let $G$ be a Borel graph and let $\mu$ be a probability \ep{Borel} measure on $V(G)$. The \emphd{$\mu$-independence number} of $G$ is the quantity $\alpha_\mu(G) \defeq \sup_I \mu(I)$, where the supremum is taken over all $\mu$-measurable $G$-independent subsets $I \subseteq V(G)$. Note that if $I_1$, \ldots, $I_\ell$ is a Borel $k$-fold $\ell$-coloring of $G$, 
	then \[\ell \alpha_\mu(G) \,\geq\, \mu(I_1) + \cdots + \mu(I_\ell) \,\geq\, k,\] which implies $\chi^\ast_\mathrm{B}(G) \geq 1/\alpha_\mu(G)$. 
	Our main result is a matching upper bound for Schreier graphs:
	
	
	
	
	\begin{theo}\label{theo:main}
        Let $\G$ be a countable group and let $F \subseteq \G$ be a finite set. If $\G \acts X$ is a free Borel action on a standard Borel space, then
        \begin{equation}\label{eq:all}
            \chi^\ast_\mathrm{B}(G(X,F)) \,\leq\, \frac{1}{\alpha_\beta(G(\Free(2^{\G}), F))}.
        \end{equation}
        In particular,
        \begin{equation}\label{eq:shift}
            \chi^\ast_\mathrm{B}(G(\Free(2^{\G}), F)) \,=\, \frac{1}{\alpha_\beta(G(\Free(2^{\G}), F))}.
        \end{equation}
	\end{theo}
	
	While \eqref{eq:shift} is a special case of \eqref{eq:all}, it is possible to deduce \eqref{eq:all} from \eqref{eq:shift} using a theorem of Seward and Tucker-Drob \cite{STD}, which asserts that every free Borel action of $\G$ admits a Borel $\G$-equivariant map to $\Free(2^\G)$. Nevertheless, we will give a simple direct proof of \eqref{eq:all} in \S\ref{sec:proof}.
	
	An interesting feature of Theorem~\ref{theo:main} is that it establishes a precise relationship between a \emph{Borel} parameter $\chi^\ast_\mathrm{B}$ and a \emph{measurable} parameter $\alpha_\beta$. We find this somewhat surprising, since ignoring sets of measure $0$ usually significantly reduces the difficulty of problems in Borel combinatorics. For instance, given a Borel graph $G$ and a probability measure $\mu$ on $V(G)$, one can consider the \emphd{$\mu$-measurable chromatic number} $\chi_\mu(G)$, i.e., the least $\ell \in \N$ such that there exist $\mu$-measurable $G$-independent sets $I_1$, \ldots, $I_\ell$ whose union is $V(G)$. By definition, $\chi_\mu(G) \leq \chi_\mathrm{B}(G)$, and it is often the case that this inequality is strict---see \cite[\S6]{KechrisMarks} for a number of examples. By contrast, as an immediate consequence of Theorem~\ref{theo:main} we obtain the opposite inequality $\chi_\mathrm{B}^\ast(G) \leq \chi_\beta(G)$, where $G$ is the Schreier graph of the free part of the shift: 
	
	\begin{corl}
	    Let $\G$ be a countable group and let $F \subseteq \G$ be a finite set. Set $G \defeq G(\Free(2^{\G}), F)$. Then $\chi^\ast_\mathrm{B}(G) \leq \chi_\beta(G)$.
	\end{corl}
	\begin{scproof}
	    The result follows from Theorem~\ref{theo:main} and the obvious inequality $\alpha_\beta(G) \geq 1/\chi_\beta(G)$. 
	\end{scproof}
	
	As a concrete application of Theorem~\ref{theo:main}, consider the free group case. For $n \geq 1$, let $\F_n$ be the free group of rank $n$ generated freely by elements $\sigma_1$, \ldots, $\sigma_n$ and let $G_n$ denote the Schreier graph of the free part of the shift action $\F_n \acts 2^{\F_n}$ with respect to the set $\set{\sigma_1, \ldots, \sigma_n}$. Then every connected component of $G_n$ is an \ep{infinite} $2n$-regular tree. In particular, the chromatic number of $G_n$ is $2$. On the other hand, Marks \cite{Marks} proved that $\chi_{\mathrm{B}}(G_n) = 2n+1$. Meehan inquired where 
	between these two extremes the Borel fractional chromatic number of $G_n$ lies: 
	
	\begin{ques}[{\cite[Question 4.6.3]{Mee}; see also \cite[Problem 8.17]{KechrisMarks}}]\label{ques:Mee}
	    What is the Borel fractional chromatic number of $G_n$? Is it always equal to $2$?
	\end{ques}
	
	Using Theorem~\ref{theo:main} together with some known results we asymptotically determine $\chi_\mathrm{B}^\ast(G_n)$ \ep{and, in particular, give a negative answer to the second part of Question~\ref{ques:Mee}}:
	
	\begin{corl}\label{corl:free}
	    For all $n \geq 1$, we have
	    \[
	        \chi_\mathrm{B}^\ast(G_n) \,=\, (2+o(1))\frac{n}{\log n},
	    \]
	    where $o(1)$ denotes a function of $n$ that approaches $0$ as $n \to \infty$.
	\end{corl}
	
	In other words, the Borel fractional chromatic number of $G_n$ is less than its ordinary Borel chromatic number roughly by a factor of $\log n$. We present the derivation of Corollary~\ref{corl:free} in \S\ref{sec:corl}.
	
	
	\section{Proof of Theorem~\ref{theo:main}}\label{sec:proof}
	
	We shall use the following theorem of Kechris, Solecki, and Todorcevic:
	
	\begin{theo}[{Kechris--Solecki--Todorcevic \cite[Proposition 4.6]{KST}}]\label{theo:KST}
	     If $G$ is a Borel graph of finite maximum degree $d$, then $\chi_\mathrm{B}(G) \leq d+1$.
	\end{theo}
	
	Fix a countable group $\G$ and a finite subset $F \subseteq \G$. Without loss of generality, we may assume that $\mathbf{1} \not \in F$. 
	%
	Say that a set $I \subseteq 2^\G$ is \emphd{independent} if $I \cap (\sigma \cdot I) = \0$ for all $\sigma \in F$ \ep{when $I \subseteq \Free(2^\G)$, this exactly means that $I$ is $G(\Free(2^\G), F)$-independent}. For brevity, let \[\alpha_\beta \,\defeq\, \alpha_\beta(G(\Free(2^\G), F)).\]
	
	\begin{lemma}\label{lemma:clopen}
	    For every $\alpha < \alpha_\beta$, there is a clopen independent set $I \subseteq 2^\G$ such that $\beta(I) \geq \alpha$.
	\end{lemma}
	\begin{scproof}
	    Let $J \subseteq \Free(2^\G)$ be a $\beta$-measurable independent set with $\beta(J) > \alpha$. Since $\beta$ is regular \cite[Theorem 17.10]{KechrisDST} and $2^\G$ is zero-dimensional, there is a clopen set $C \subseteq 2^\G$ with
	    \[\mu(J \symdif C) \,\leq\, \frac{\beta(J) - \alpha}{|F| + 1}.\] Set $I \defeq C \setminus \bigcup_{\sigma \in F} (\sigma \cdot C)$. By construction, $I$ is clopen and independent. Moreover, if $x \in J \setminus I$, then either $x \in J \setminus C$ or $x \in (\sigma \cdot C) \setminus (\sigma \cdot J)$ for some $\sigma \in F$. Therefore, \[\beta(I) \,\geq\, \beta(J) \,-\, (|F|+1) \beta(J \symdif C) \,\geq\, \alpha. \qedhere\]
	\end{scproof}
	
	Let $\G \acts X$ be a free Borel action on a standard Borel space. Fix an arbitrary clopen independent set $I \subseteq 2^\G$. We will prove that $\chi^\ast_\mathrm{B}(G(X,F)) \leq 1/\beta(I)$, which yields Theorem~\ref{theo:main} by Lemma~\ref{lemma:clopen}. Since $I$ is clopen, there exist finite sets $D \subseteq \G$ and $\Phi \subseteq 2^D$ such that
	\[
	    I \,=\, \set{x \in 2^\G \,:\, \rest{x}{D} \in \Phi},
	\]
	where $\rest{x}{D}$ denotes the restriction of $x$ to $D$. Note that
	\[
	    \beta(I) \,=\, \frac{|\Phi|}{2^{|D|}}.
	\]
	Let $N \defeq |DD^{-1}|$ and consider the graph $H \defeq G(X, DD^{-1})$. Every vertex in $H$ has precisely $N-1$ neighbors \ep{we are subtracting $1$ to account for the fact that a vertex is not adjacent to itself}. By Theorem~\ref{theo:KST}, this implies that $\chi_\mathrm{B}(H) \leq N$. 
	In other words, we may fix a Borel function $f \colon X \to N$ such that $f(u) \neq f(v)$ whenever $u$, $v \in X$ are distinct points satisfying $v \in DD^{-1} \cdot u$. This implies that for each $x \in X$, 
	the restriction of $f$ to the set $D \cdot x$ is injective. 
	Now, to each mapping $\phi \colon N \to 2$, we associate a Borel $\G$-equivariant function $\pi_\phi \colon X \to 2^\G$ as follows:
	\[
	    \pi_\phi(x)(\gamma) \,\defeq\, (\phi \circ f)(\gamma \cdot x) \quad \text{for all } x \in X \text{ and } \gamma \in \G.
	\]
	Let $I_\phi \defeq \pi_\phi^{-1}(I)$. Since $\pi_\phi$ is $\G$-equivariant, $I_\phi$ is $G(X, F)$-independent. Consider any $x \in X$ and let
	\[
	    S_x \,\defeq\, \set{f(\gamma \cdot x) \,:\, \gamma \in D}.
	\]
	By the choice of $f$, $S_x$ is a subset of $N$ of size $|D|$. Whether or not $x$ is in $I_\phi$ is determined by the restriction of $\phi$ to $S_x$; furthermore, there are exactly $|\Phi|$ such restrictions that put $x$ in $I_\phi$. Thus, the number of mappings $\phi \colon N \to 2$ for which $x \in I_\phi$ is \[|\Phi|2^{N - |D|} \,=\, \beta(I)2^N.\] Since this holds for all $x \in X$, we conclude that the sets $I_\phi$ cover every point in $X$ exactly $\beta(I)2^N$ times. Therefore, $\chi^\ast_\mathrm{B}(G(X,F)) \leq 1/\beta(I)$, as desired.
	
	\section{Proof of Corollary~\ref{corl:free}}\label{sec:corl}
	
	  Thanks to Theorem~\ref{theo:main}, in order to establish Corollary~\ref{corl:free} it is enough to verify that
	   \[
	        \alpha_\beta(G_n) \,=\, \left(\frac{1}{2} + o(1)\right) \frac{\log n}{n}.
	   \]
	   There are a number of known constructions that witness the lower bound \[\alpha_\beta(G_n) \,\geq\, \left(\frac{1}{2} + o(1)\right) \frac{\log n}{n};\] see, e.g.,  \cite{LW} by Lauer and Wormald and \cite{GG} by Gamarnik and Goldberg. Moreover, by \cite[Corollary 1.2]{Ber}, even the inequality $\chi_\beta(G_n) \leq (2+o(1))n/\log n$ holds. For the upper bound
	   \begin{equation}\label{eq:upper}
	    \alpha_\beta(G_n) \,\leq\, \left(\frac{1}{2} + o(1)\right)\frac{\log n}{n},
	   \end{equation}
	   we shall use a theorem of Rahman and Vir\'{a}g \cite{RV}, which says that the largest density of a factor of i.i.d.~independent set in the $d$-regular tree is at most $(1+o(1))\log d/d$. In the remainder of this section we describe their result and explain how it implies the desired upper bound on $\alpha_\beta(G_n)$. 
	    
	   Fix an integer $n \geq 1$. For our purposes, it will be somewhat more convenient to work on the space $[0,1]^{\F_n}$ instead of $2^{\F_n}$, where $[0,1]$ is the unit interval equipped with the usual Lebesgue probability measure. The product measure on $[0,1]^{\F_n}$ is denoted by $\lambda$. Let $H_n$ be the Schreier graph of the shift action $\F_n \acts [0,1]^{\F_n}$ corresponding to the standard generating set of $\F_n$. We remark that, by a theorem of Ab\'ert and Weiss \cite{AW} \ep{see also \cite[Theorem 6.46]{KechrisMarks}}, $\alpha_\beta(G_n) = \alpha_\lambda(H_n)$, so it does not really matter whether we are working with $G_n$ or $H_n$. 
	   
	   Set $d \defeq 2n$ and let $\mathbb{T}_d$ denote the Cayley graph of the free group $\F_n$ with respect to the standard generating set. In other words, $\mathbb{T}_d$ is an \ep{infinite} $d$-regular tree. We view $\mathbb{T}_d$ as a \emph{rooted} tree, whose root is the vertex $\mathbf{1}$, i.e., the identity element of $\F_n$. Let $\mathfrak{A}$ be the automorphism group of $\mathbb{T}_d$, i.e., the set of all bijections $\mathfrak{a} \colon \F_n \to \F_n$ that preserve the edges of $\mathbb{T}_d$, and let $\mathfrak{A}_\bullet \subseteq \mathfrak{A}$ be the subgroup comprising the root-preserving automorphisms, i.e., those $\mathfrak{a} \in \mathfrak{A}$ that map $\mathbf{1}$ to $\mathbf{1}$. The space $[0,1]^{\F_n}$ is equipped with a natural right action $[0,1]^{\F_n} \racts \mathfrak{A}$. Namely, for $\mathfrak{a} \in \mathfrak{A}$ and $x \in [0,1]^{\F_n}$, the result of acting by $\mathfrak{a}$ on $x$ is the function $x \cdot \mathfrak{a} \colon \F_n \to [0,1]$ given by
	   \[
	        (x \cdot \mathfrak{a})(\delta) \,\defeq\, x(\mathfrak{a}(\delta)) \quad \text{for all } \delta \in \F_n.
	   \]
	   For each $\gamma \in \F_n$, there is a corresponding automorphism $\mathfrak{a}_\gamma \in \mathfrak{A}$ sending every group element $\delta \in \F_n$ to $\delta \gamma$. The mapping $\F_n \to \mathfrak{A} \colon \gamma \mapsto \mathfrak{a}_\gamma$ is an antihomomorphism of groups, that is, we have \[\mathfrak{a}_{\gamma \sigma} \,=\, \mathfrak{a}_\sigma \circ \mathfrak{a}_\gamma \quad \text{for all $\gamma$, $\sigma \in \F_n$},\] where $\circ$ denotes composition. In particular, $\set{\mathfrak{a}_\gamma \,:\, \gamma \in \F_n}$ is a subgroup of $\mathfrak{A}$ isomorphic to $\F_n$. The right action $[0,1]^{\F_n} \racts \mathfrak{A}$ and the left action $\F_n \acts [0,1]^{\F_n}$ are related by the formula
	   \[
	    x \cdot \mathfrak{a}_\gamma = \gamma \cdot x \quad \text{for all $x \in [0,1]^{\F_n}$}.
	   \] 
	   A set $X \subseteq [0,1]^{\F_n}$ is called \emphd{$\mathfrak{A}_\bullet$-invariant} if 
	   $x \cdot \mathfrak{a} \in X$ for all $x \in X$ and $\mathfrak{a} \in \mathfrak{A}_\bullet$. The Rahman--Vir\'{a}g theorem can now be stated as follows:
	   
	   \begin{theo}[{Rahman--Vir\'{a}g \cite[Theorem 2.1]{RV}}]\label{theo:RV}
	        If $I \subseteq [0,1]^{\F_n}$ is an $\mathfrak{A}_\bullet$-invariant $\lambda$-measurable $H_n$-independent set, then
	        \[
	            \lambda(I) \,\leq\, (1+o(1)) \frac{\log d}{d} \,=\, \left(\frac{1}{2} + o(1)\right) \frac{\log n}{n}.
	        \]
	   \end{theo}
	   
	   Theorem~\ref{theo:RV} is almost the result we want, except that we need an upper bound on the measure of \emph{every} \ep{not necessarily $\mathfrak{A}_\bullet$-invariant} $\lambda$-measurable $H_n$-independent set $I$. To remove the $\mathfrak{A}_\bullet$-invariance assumption, we use the following consequence of Theorem~\ref{theo:RV}:
	   
	   \begin{corl}\label{corl:tree}
	       There exists a Borel graph $Q$ with a probability measure $\mu$ on $V(Q)$ such that:
	       \begin{itemize}
	           \item every connected component of $Q$ is a $d$-regular tree; and
	           \item $\alpha_\mu(Q) \leq (1/2 + o(1)) \log n/n$.
	       \end{itemize}
	   \end{corl}
	   \begin{scproof}
	    We use a construction that was studied by Conley, Kechris, and Tucker-Drob in \cite{Ultra}. Let $\Omega$ be the set of all points $x \in [0,1]^{\F_n}$ such that $x \cdot \mathfrak{a} \neq x$ for every non-identity automorphism $\mathfrak{a} \in \mathfrak{A}$. Let us make a couple observations about $\Omega$. Notice that, by definition, the set $\Omega$ is invariant under the action $[0,1]^{\F_n} \racts \mathfrak{A}$; in particular, it is invariant under the shift action $\F_n \acts [0,1]^{\F_n}$. Furthermore, the induced action of $\F_n$ on $\Omega$ is free \ep{indeed, even the action $\Omega \racts \mathfrak{A}$ is free}. Since every injective mapping $\F_n \to [0,1]$ belongs to $\Omega$, we conclude that $\lambda(\Omega) = 1$.
	    Now consider the quotient space $V \defeq \Omega/\mathfrak{A}_\bullet$. As the group $\mathfrak{A}_\bullet$ is compact, the space $V$ is standard Borel \cite[paragraph preceding Lemma 7.8]{Ultra}. Let $\mu$ be the push-forward of $\lambda$ under the quotient map $\Omega \to V$, and let $Q$ be the graph with vertex set $V$ in which two vertices $\bm{x}$, $\bm{y} \in V$ are adjacent if and only if there are representatives $x \in \bm{x}$ and $y \in \bm{y}$ that are adjacent in $H_n$. Conley, Kechris, and Tucker-Drob \cite[Lemma 7.9]{Ultra} \ep{see also \cite[Proposition 1.9]{Thornton}} showed that every connected component of $Q$ is a $d$-regular tree. Furthermore, by construction, a set $I \subseteq V$ is $Q$-independent if and only if its preimage under the quotient map is $H_n$-independent. Since the quotient map establishes a one-to-one correspondence between subsets of $V$ and $\mathfrak{A}_\bullet$-invariant subsets of $\Omega$, Theorem~\ref{theo:RV} is equivalent to the assertion that $\alpha_\mu(Q) \leq (1/2 + o(1)) \log n/n$, as desired.
	   \end{scproof}
	   
	   In view of Corollary~\ref{corl:tree}, the following lemma completes the proof of \eqref{eq:upper}:
	   
	   \begin{lemma}\label{lemma:AW}
	        Let $Q$ be a Borel graph in which every connected component is a $d$-regular tree and let $\mu$ be a probability measure on $V(Q)$. Then $\alpha_\mu(Q) \geq \alpha_\beta(G_n)$.
	   \end{lemma}
	   
	   In the case when $Q$ is the Schreier graph of a free measure-preserving action of $\F_n$, the conclusion of Lemma~\ref{lemma:AW} follows from the Ab\'ert--Weiss theorem~\cite{AW}. To handle the general case, we rely on a strengthening of a recent result of T\'oth \cite{Toth} due to Greb\'ik \cite{Jan}, which, roughly, asserts that every $d$-regular Borel graph is ``approximately'' induced by an action of $\F_n$. 

	   To state this result precisely, we introduce the following terminology. 
	   A \emphd{Borel partial action} $\bm{p}$ of $\F_n$ on a standard Borel space $X$, in symbols $\bm{p} \colon \F_n \acts^\ast X$, is a sequence of Borel partial injections $p_1$, \ldots, $p_n \colon X \pto X$. 
	   Given a Borel graph $Q$, we say that a Borel partial action $\bm{p} \colon \F_n \acts^\ast V(Q)$ is a \emphd{partial Schreier decoration} of $Q$ if $p_i(x)$ is adjacent to $x$ for all $1 \leq i \leq n$ and $x \in \dom(p_i)$. If $\bm{p}$ is a partial Schreier decoration of a graph $Q$, then we let $C(Q, \bm{p})$ be the set of all vertices $x \in V(Q)$ such that $x$ belongs to both the domain and the image of every $p_i$ and the neighborhood of $x$ in $Q$ is equal to the set $\set{p_1(x), \ldots, p_n(x), p^{-1}_1(x), \ldots, p_n^{-1}(x)}$. A \emphd{Schreier decoration} of $Q$ is a partial Schreier decoration $\bm{p}$ such that $C(Q, \bm{p}) = V(Q)$. It is easy to see that $Q$ admits a Schreier decoration if and only if it is the Schreier graph of a Borel action of $\F_n$. 
	   
	   Now we can state Greb\'ik's result: 
	   
	   \begin{theo}[{Greb\'ik \cite[Theorem 0.2(III)]{Jan}}]\label{theo:Jan}
	        Let $Q$ be a $d$-regular Borel graph and let $\mu$ be a probability measure on $V(Q)$. Then for every $\epsilon > 0$, $Q$ admits a partial Schreier decoration $\bm{p}$ such that $\mu(C(Q, \bm{p})) \geq 1 - \epsilon$.
	   \end{theo}
	   
	   With Theorem~\ref{theo:Jan} in hand, we are ready to establish Lemma~\ref{lemma:AW}.
	   
	   \begin{scproof}[ of Lemma~\ref{lemma:AW}]
	        Recall that we denote the generators of $\F_n$ by $\sigma_1$, \ldots, $\sigma_n$. Let $Q$ be a Borel graph in which every connected component is a $d$-regular tree and let $\mu$ be a probability measure on $V(Q)$. Thanks to Lemma~\ref{lemma:clopen}, it suffices to show that $\alpha_\mu(Q) \geq \beta(I)$ for every clopen independent set $I \subseteq 2^{\F_n}$, where, as in \S\ref{sec:proof}, we say that $I$ is independent if $I \cap (\sigma_i \cdot I) = \0$ for each $1 \leq i \leq n$.
	        
	        Fix a clopen independent set $I \subseteq 2^{\F_n}$. Since $I$ is clopen, we can write
	        \[
	            I \,=\, \set{x \in 2^{\F_n} \,:\, \rest{x}{D} \in \Phi},
	        \]
	        where $D \subset \F_n$ and $\Phi \subseteq 2^D$ are finite sets. Furthermore, we may assume without loss of generality that $D = \set{\gamma \in \F_n \,:\, |\gamma| \leq k}$ for some $k \in \N$, where $|\gamma|$ denotes the word norm of $\gamma$. For a vertex $x \in V(Q)$, we let $N^k(x)$ be the set of all vertices that are joined to $x$ by a path of length at most $k$. Since every connected component of $Q$ is a $d$-regular tree, we have $|N^k(x)| = |D|$ for all $x \in V(Q)$. This allows us to define a probability measure $\mu_k$ on $V(Q)$ via
	        \[
	            \mu_k(A) \,\defeq\, \int \frac{|A \cap N^k(x)|}{|D|} \,\mathrm{d}\mu(x) \quad \text{for all Borel } A \subseteq V(Q). 
	        \]
	       
	       We have now prepared the ground for an application of Theorem~\ref{theo:Jan}. Fix $\epsilon > 0$ and let $\bm{p}$ be a partial Schreier decoration of $Q$ such that
	       \[
	            \mu_k(C(Q, \bm{p})) \,\geq\, 1 - \frac{\epsilon}{|D|},
	       \]
	       which exists by Theorem~\ref{theo:Jan}. Let $C_k$ be the set of all $x \in V(Q)$ such that $N^k(x) \subseteq C(Q, \bm{p})$. Then
	       \begin{align*}
	           1 - \frac{\epsilon}{|D|} \,\leq\, \mu_k(C(Q, \bm{p})) \,&=\, \int \frac{|C(Q, \bm{p}) \cap N^k(x)|}{|D|} \,\mathrm{d}\mu(x) \\
	           &\leq\, \mu(C_k) + \left(1 - \frac{1}{|D|}\right)(1 - \mu(C_k))\,=\, \frac{1}{|D|} \mu(C_k) + 1 - \frac{1}{|D|},
	       \end{align*}
	       which implies that $\mu(C_k) \geq 1 - \epsilon$. The importance of the set $C_k$ lies in the fact that for each $x \in C_k$ and $\gamma \in D$, there is a natural way to define the notation $\gamma \cdot x$. Namely, we write $\gamma$ as a reduced word: \[\gamma \,=\, \sigma_{i_1}^{s_1} \cdots \sigma_{i_\ell}^{s_\ell},\] where $0 \leq \ell \leq k$, each index $i_j$ is between $1$ and $n$, and each $s_j$ is $1$ or $-1$. Since $N^k(x) \subseteq C(Q, \bm{p})$, there is a unique sequence $x_0$, $x_1$, \ldots, $x_\ell$ of vertices with
	       \[
	        x_0 = x \quad \text{and} \quad x_{j} = p_{i_j}^{s_j}(x_{j-1}) \text{ for all } 1 \leq j \leq \ell.
	       \]
	       We then set $\gamma \cdot x \defeq x_\ell$. Note that we have $N^k(x) = \set{\gamma \cdot x \,:\, \gamma \in D}$.
	       
	       The remainder of the argument utilizes a construction similar to the one in the proof of Theorem~\ref{theo:main} given in \S\ref{sec:proof}. Consider the graph $R$ with the same vertex set as $Q$ in which two distinct vertices are adjacent if and only if they are joined by a path of length at most $2k$ in $Q$. Since every connected component of $Q$ is a $d$-regular tree, each vertex in $R$ has the same finite number of neighbors, so, by Theorem~\ref{theo:KST}, the Borel chromatic number $\chi_\mathrm{B}(R)$ is finite. Let $N \defeq \chi_\mathrm{B}(R)$ and fix a Borel function $f \colon V(Q) \to N$ such that $f(u) \neq f(v)$ whenever $u$ and $v$ are adjacent in $R$. Then for each $x \in V(Q)$, 
	        the restriction of $f$ to the set $N^k(x)$ is injective. 
	        Now, to each mapping $\phi \colon N \to 2$, we associate function $\pi_\phi \colon C_k \to 2^D$ as follows:
	        \[
	            \pi_\phi(x)(\gamma) \,\defeq\, (\phi \circ f)(\gamma \cdot x) \quad \text{for all } x \in C_k \text{ and } \gamma \in D.
	        \]
	        Let $I_\phi \defeq \set{x \in C_k \,:\, \pi_\phi(x) \in \Phi}$. The independence of $I$ implies that the set $I_\phi$ is $Q$-independent. We will show that for some choice of $\phi \colon N \to 2$, $\mu(I_\phi) \geq (1-\epsilon) \beta(I)$. Since $\epsilon$ is arbitrary, this yields the desired bound $\alpha_\mu(Q) \geq \beta(I)$ and completes the proof of Lemma~\ref{lemma:AW}.
	        
	        Consider any $x \in C_k$ and let
        	\[
        	    S_x \,\defeq\, \set{f(\gamma \cdot x) \,:\, \gamma \in D}.
        	\]
        	Since $f$ is injective on $N^k(x)$, $S_x$ is a subset of $N$ of size $|D|$. Whether or not $x$ is in $I_\phi$ is determined by the restriction of $\phi$ to $S_x$; furthermore, there are exactly $|\Phi|$ such restrictions that put $x$ in $I_\phi$. Thus, the number of mappings $\phi \colon N \to 2$ for which $x \in I_\phi$ is \[|\Phi|2^{N - |D|} \,=\, \beta(I)2^N.\] Since this holds for all $x \in C_k$, we conclude that
        	\[
        	    \sum_{\phi \colon N \to 2} \mu(I_\phi) \,\geq\, \mu(C_k) \beta(I)2^N \,\geq\, (1-\epsilon)\beta(I)2^N,
        	\]
        	where the second inequality uses that $\mu(C_k) \geq 1 - \epsilon$. In other words, the average value of $\mu(I_\phi)$ over all $\phi \colon N \to 2$ is at least $(1-\epsilon) \beta(I)$. Thus, the maximum is at least $(1-\epsilon) \beta(I)$ as well, and the proof is complete.
	   \end{scproof}
	   
	   \subsubsection*{Acknowledgement} I am grateful to the anonymous referees for carefully reading this paper and providing helpful feedback.

	
	\printbibliography

@book{FracBook,
	author = {E.R. Scheinerman and D.H. Ullman},
	title = {Fractional Graph Theory},
	date = {1997},
	publisher = {John Wiley \& Sons},
	url = {https://www.ams.jhu.edu/ers/books/fractional-graph-theory-a-rational-approach-to-the-theory-of-graphs/},
}

@unpublished{KechrisMarks,
	author = {A.S. Kechris and A.S. Marks},
	title = {Descriptive Graph Combinatorics},
	date = {2020},
	howpublished = {\url{http://www.math.caltech.edu/~kechris/papers/combinatorics20book.pdf} (preprint)},
}

@unpublished{Pikh_survey,
	author = {O. Pikhurko},
	title = {Borel combinatorics of locally finite graphs},
	howpublished = {\url{https://arxiv.org/abs/2009.09113} (preprint)},
	date = {2020},
}

@thesis{Mee,
	author = {C. Meehan},
	title = {Definable combinatorics of graphs and equivalence relations},
	type = {Ph.D. Thesis},
	institution = {California Institute of Technology},
	date = {2018},
	location = {Pasadena, CA},
	url = {https://resolver.caltech.edu/CaltechTHESIS:06012018-160828760},
}

@article{LW,
	author = {J. Lauer and N. Wormald},
	title = {Large independent sets in regular graphs of large girth},
	journaltitle = {J. Combin. Theory},
	series = {B},
	volume = {97},
	pages = {999--1009},
	date = {2007},
}

@article{RV,
	author = {M. Rahman and B. Vir{\'{a}}g},
	title = {Local algorithms for independent sets are half-optimal},
	journaltitle = {Ann. Probab.},
	volume = {45},
	number = {3},
	pages = {1543--1577},
	date = {2017},
}

@unpublished{Thornton,
	author = {R. Thornton},
	title = {Factor of i.i.d. processes and Cayley diagrams},
	howpublished = {\url{https://arxiv.org/abs/2011.14604} (preprint)},
	date = {2021},
}

@article{AW,
	author = {M. Ab{\'{e}}rt and B. Weiss},
	title = {Bernoulli actions are weakly contained in any free action},
	journaltitle = {Ergod. Th. and Dynam. Sys.},
	volume = {33},
	number = {2},
	date = {2013},
	page = {323--333},
}

@article{Ultra,
	author = {C.T. Conley and A.S. Kechris and R.D. Tucker-Drob},
	title = {Ultraproducts of measure preserving actions and graph combinatorics},
	journaltitle = {Ergod. Th. and Dynam. Sys.},
	volume = {33},
	number = {2},
	date = {2013},
	page = {334--374},
}

@article{Toth,
	author = {L.M. T{\'{o}}th},
	title = {Invariant Schreier decorations of unimodular random networks},
	date = {2021},
	volume = {4},
	journaltitle = {Ann. H. Lebesgue},
	pages = {1705--1726},
}

@unpublished{Jan,
	author = {J. Greb{\'{i}}k},
	title = {Approximate Schreier decorations and approximate K{\"{o}}nig’s line coloring Theorem},
	howpublished = {\url{https://users.math.cas.cz/~grebik/Approx.pdf} (preprint)},
	date = {2020},
	addendum = {Ann. H. Lebesgue (to appear)},
}

@article{GG,
	author = {D. Gamarnik and D.A. Goldberg},
	title = {Randomized greedy algorithms for independent sets and matchings in regular graphs: Exact results and finite girth corrections},
	journaltitle = {Combin. Probab. Comput.},
	volume = {19},
	date = {2010},
	page = {61--85},
}

@article{Ber,
	author = {A. Bernshteyn},
	title = {Measurable versions of the Lov{\'{a}}sz Local Lemma and measurable graph colorings},
	journaltitle = {Adv. Math.},
	volume = {353},
	date = {2019},
	page = {153--223},
}

@article{STD,
	author = {B. Seward and R.D. Tucker-Drob},
	title = {Borel structurability on the $2$-shift of a countable group},
	journaltitle = {Ann. Pure Appl. Logic},
	volume = {167},
	number = {1},
	date = {2016},
	pages = {1--21},
}

@book{KechrisDST,
	author = {A.S. Kechris},
	title = {Classical Descriptive Set Theory},
	date = {1995},
	publisher = {Springer-Verlag},
	location = {New York},
}

@article{KST,
	author = {A.S. Kechris and S. Solecki and S. Todorcevic},
	title = {Borel chromatic numbers},
	journaltitle = {Adv. in Math.},
	date = {1999},
	volume = {141},
	pages = {1--44},
}

@article{Marks,
	author = {A.S. Marks},
	title = {A determinacy approach to Borel combinatorics},
	journaltitle = {J. Amer. Math. Soc.},
	volume = {29},
	date = {2016},
	pages = {579--600},
}
    
\end{document}